\documentclass[12pt,a4paper]{article}
\usepackage{color}
\usepackage{authblk}
\usepackage{cmap}
\usepackage{amsmath,amssymb}
\usepackage[T1]{fontenc}
\usepackage[utf8]{inputenc}
\usepackage[english]{babel}
\usepackage{hyperref}

\begin{document}
\author[1,2,*]{V.\,G.~Zborovskii}
\affil[1]{Troitsk Branch of the P.\,N.~Lebedev Physical Institute of the Russian Academy of Sciences, Troitsk, Moscow, Russia}
\affil[2]{National Research Centre ``Kurchatov Institute'', Moscow, Russia}
\affil[*]{E-mail: \texttt{ZborovskiyVG@lebedev.ru}}
\title{On Some Special Cases of the Weber--Schafheitlin Integral}
\date{2026}
\maketitle

\begin{abstract}
We consider special cases of the discontinuous Weber--Schafheitlin integral
in which the orders $\mu$ and $\nu$ of the Bessel functions differ by an integer
and the exponent $\lambda$ equals $0$ or $\pm1$, with $\mu-\nu+\lambda$ odd.
Explicit closed-form results are obtained in terms of Jacobi polynomials.

\end{abstract}

\section{Introduction}
\label{sec:intro}

In~\cite{Zborovskii25}, a spectral method of spherical harmonics in the Fourier
representation was proposed for computing the elastic fields, namely, displacements,
strains, stresses, and body forces, produced by ellipsoidal inclusions in
anisotropic media. The method employs integrals of the form
\begin{equation}
	\label{eq:HankelRBF}
	R_{\lambda,k,l} (z)
	= z^{-\frac12} \int_0^\infty
	\xi^{-\lambda}J_{k+\frac32}(\xi) J_{l+\frac12}(\xi z)\, d\xi,
\end{equation}
which are Hankel transforms of half-integer order applied to radial basis functions.
The boundary of the inclusion is given by $z=1$. In~\eqref{eq:HankelRBF}
the indices $k$ and $l$ are integers; $\lambda=0$ for strains and stresses,
$\lambda=1$ for displacements, and $\lambda=-1$ for body forces per unit volume.
For brevity, the indices $k$ and $l$ will sometimes be omitted.
For the spectral method of~\cite{Zborovskii25} it suffices to consider the case
\begin{equation}
	\label{eq:spectral_case}
	l = k + 2p + s,\quad
	s = |\lambda|,\quad
	p \in \mathbb{Z},\quad
	k \ge -1,\quad
	l \ge 0.
\end{equation}

For $\lambda>-1$, the integral~\eqref{eq:HankelRBF} belongs to the class of
discontinuous Weber--Schafheitlin integrals. The general case was studied in
Watson's monograph~\cite[Secs.~13.4, 13.41]{Watson}, where it is shown that
the values of the integral are expressed in terms of the Gauss hypergeometric
function.

We study the integral~\eqref{eq:HankelRBF} as a function of $z$ for fixed
$k$ and $l$ and for various values of $\lambda$.
The integral~\eqref{eq:HankelRBF} converges at $\xi \to 0$ provided
$k+l+2-\lambda \ge 0$. Under the constraints~\eqref{eq:spectral_case}
this condition is always satisfied.

The properties of the integral~\eqref{eq:HankelRBF} are determined by its
convergence at the upper limit $\xi \to \infty$. For $\lambda>0$, the
integral~\eqref{eq:HankelRBF} converges absolutely for all $z$, including
$z=1$, and $R_\lambda(z)$ is continuous~\cite{Watson}.

For $\lambda=0$ and $l=k+2p$, $p \in \mathbb{Z}$, the
integral~\eqref{eq:HankelRBF} converges (conditionally) for all $z$,
including $z=1$; however, the limits of $R_\lambda(z)$ as $z \to 1\pm$ do not
coincide, and the value $R_\lambda(1)$ equals the half-sum of these
limits~\cite{Watson}.

Now consider the integral~\eqref{eq:HankelRBF} for $s=1$ and $\lambda \le 0$.
For $z \ne 1$ the integral converges (conditionally) for $\lambda > -1$ and
diverges for $\lambda=-1$; however, as will be shown below, the limit of
$R_\lambda(z)$ as $\lambda \to -1$ exists. At the same time, $R_\lambda(1)$
diverges for all $\lambda < 0$. It is therefore natural to define
$R_{-1}(z) = \lim_{\lambda \to -1} R_\lambda(z)$ and to treat this limit as a
distribution whose singular part is supported on the set $z=1$.
The limit $R_{-1}(z)$ was studied in~\cite{Maximon,KR}, where the explicit form
of the singular term was obtained. For a special case
of the form~\eqref{eq:HankelRBF} under conditions~\eqref{eq:spectral_case},
the regular part of
the integral~\eqref{eq:HankelRBF} was expressed in~\cite{Maximon} in terms of
derivatives of Legendre polynomials.

In the present paper, simple closed-form results are obtained for the special
cases of the Weber--Schafheitlin integral under discussion, with the regular
part expressed in terms of Jacobi polynomials.
%, and the singular part in terms of the Dirac delta function $\delta(z-1)$.
Compared to higher-order derivatives of Legendre polynomials~\cite{Maximon},
 these expressions are convenient for computation 
 whether used individually or in series expansions~\cite{NR}.

\section{General form of the integral}

For given $k$, $p$, $s$ and $\lambda>-1$, the integral~\eqref{eq:HankelRBF}
is expressed in terms of the Gauss hypergeometric function $F$~\cite{Watson}.

Denote by $R^<_{\lambda,k,l} (z)$ the value of the integral $R$ for $z < 1$. Then
\begin{equation}
	\label{eq:hypergeo_int}
	R^<_{\lambda,k,l} (z) = \frac{z^{\gamma -\frac32}}
	{2^{\lambda}}
	\frac{\Gamma(\alpha)}{\Gamma(1-\beta)\Gamma(\gamma)}
	F\left(\alpha,\beta,\gamma;z^2\right),
\end{equation}
where
\begin{gather*}
	\alpha=k+p+\frac{s+3-\lambda}2,\quad
	\beta=p+\frac{s-\lambda}2,\quad
	\gamma=\alpha+\beta+\lambda=l+\frac32.
\end{gather*}

The value $R^>_{\lambda,k,l} (z)$ of the integral~\eqref{eq:HankelRBF} for
$z > 1$ is obtained from~\eqref{eq:HankelRBF} by the change of variables
$\eta=\xi z$, $w=1/z$:
\begin{equation}
	\label{eq:ext_int}
	R^>_{\lambda,k,l} (z) =
	z^{\lambda-2}
	R^<_{\lambda,l-1,k+1} \left(z^{-1}\right),
\end{equation}
where the formulas for $R^<_\lambda$ are now to be evaluated with
$p^> = 1 -s -p$ in place of $p$.

Combining~\eqref{eq:hypergeo_int} and~\eqref{eq:ext_int}, we obtain the
following expression for the integral $R$:
\begin{equation}
	\label{eq:reg}
	R_{\lambda,k,l}(z) = R^<_{\lambda,k,l}(z)\,\theta(1-z) +
	R^>_{\lambda,k,l}(z)\,\theta(z-1),
\end{equation}
where $\theta(x)$ is the Heaviside step function. For $\lambda=-1$,
Eq.~\eqref{eq:reg} gives only the regular part; see
Section~\ref{sec:f}.

It is important for the subsequent calculations that $\alpha$ and $\gamma$ do
not vanish for the parameter values under consideration.

\section{Evaluation of the integral for strains}
\label{sec:e}

The simplest to evaluate is the integral $R_0(z)$, which corresponds to the
case of strains and stresses. Substitution into~\eqref{eq:hypergeo_int} gives
\begin{equation*}
	R^<_0(z) = z^l
	\frac{\Gamma(\alpha)}{\Gamma(1-p)\Gamma(\gamma)}
	F\left(\alpha,p,\gamma;z^2\right).
\end{equation*}
For $p>0$, i.e., $l>k$, the whole expression vanishes owing to the pole of
$\Gamma(1-p)$. For $p \le 0$, the hypergeometric series terminates after the
term of order $-p$ and reduces to a Jacobi
polynomial~\cite[Eq.~18.5.7]{DLMF}:
\begin{equation*}
	F\left(\alpha,p,\gamma;z^2\right) =
	\frac{\Gamma(1-p)\Gamma(\gamma)}{\Gamma(\alpha)}
	P^{(l+\frac12,0)}_{-p}(1-2z^2),
\end{equation*}
which yields
\begin{equation}
	\label{eq:e_int}
	R^<_{0,k,l}(z) =
	\begin{cases}
		z^l P^{(l+\frac12,0)}_{-p}(1-2z^2), & p \le 0,\\
		0, & p > 0.
	\end{cases}
\end{equation}

For $z>1$, transformation~\eqref{eq:ext_int} gives
\begin{equation}
	\label{eq:e_ext}
	R^>_{0,k,l}(z) =
	\begin{cases}
		z^{-k-3} P^{(k+\frac32,0)}_{p-1}(1-2z^{-2}), & p > 0,\\
		0, & p \le 0.
	\end{cases}
\end{equation}

Evaluating the one-sided limits as $z \to 1\pm$ gives $R^<_0(1-) = (-1)^p$
for $p \le 0$ and $R^>_0(1+) = -(-1)^p$ for $p > 0$; for the other values of
$p$ the corresponding limits are zero. Thus, strains and stresses turn out to
be discontinuous functions of $z$. Physically, the jump in stresses (and
strains) is caused by concentrated forces distributed as a single layer on the
inclusion boundary $z=1$.

Although the constraint $l \ge 0$ was stated in Section~\ref{sec:intro},
Eqs.~\eqref{eq:e_int} and~\eqref{eq:e_ext} also remain valid for
$l=-1$, provided $k \ge -1$.

\section{Evaluation of the integral for displacements}

For the case of displacements, we evaluate the integral $R_1(z)$. Substitution
into~\eqref{eq:hypergeo_int} gives
\begin{equation*}
	R^<_1(z) = \frac{z^l}2
	\frac{\Gamma(\alpha)}{\Gamma(1-p)\Gamma(\gamma)}
	F\left(\alpha,p,\gamma;z^2\right).
\end{equation*}
For $p>0$, the whole expression vanishes owing to the pole of $\Gamma(1-p)$.
For $p=0$ ($l=k+1$), the hypergeometric series reduces to $1$, hence
\begin{equation*}
	R^<_{1,k,k+1}(z) = \frac{z^{k+1}}{2k+3}.
\end{equation*}
For $p<0$, to obtain Jacobi polynomials one uses the contiguous relation for
the hypergeometric function~\cite[Eq.~15.5.13]{DLMF} with $b$ replaced by
$b+1$. After this transformation, the hypergeometric series reduces to a
Jacobi polynomial of degree $-p-1$:
\begin{equation}
	\label{eq:u_int}
	R^<_{1,k,l}(z) =
	\begin{cases}
		\dfrac{z^l}{2(-p)}
		(1-z^2)
		P^{(l+\frac12,1)}_{-p-1}(1-2z^2), & p < 0,\\[2ex]
		\dfrac{z^{k+1}}{2k+3}, & p = 0,\\[2ex]
		0, & p > 0.
	\end{cases}
\end{equation}

For $z>1$, transformation~\eqref{eq:ext_int} gives
\begin{equation}
	\label{eq:u_ext}
	R^>_{1,k,l}(z) =
	\begin{cases}
		\dfrac{z^{-k-2}}{2p}
		(1-z^{-2})
		P^{(k+\frac32,1)}_{p-1}(1-2z^{-2}), & p > 0,\\[2ex]
		\dfrac{z^{-k-2}}{2k+3}, & p = 0,\\[2ex]
		0, & p < 0.
	\end{cases}
\end{equation}

The result given by Eq.~\eqref{eq:u_int} and~\eqref{eq:u_ext} is
continuous at $z=1$ (on the inclusion boundary), in agreement with the general
properties of the integral~\eqref{eq:HankelRBF}. The terms with $p \ne 0$ make
no contribution at the boundary, since they vanish there. The required
continuity of the displacements across the inclusion boundary is provided by
the terms with $p = 0$.

\section{Evaluation of the integral for body forces}
\label{sec:f}

The case of body forces corresponds to $\lambda \to -1+$.
We first compute the regular part $R^{\mathrm{reg}}_{-1,k,l}(z)$ for $z \ne 1$.
Here $s=1$; let $\lambda=-1+2\varepsilon$, where $\varepsilon \to 0+$.
Substitution into~\eqref{eq:hypergeo_int} gives
\begin{equation*}
	R^<_\lambda(z) = 2^{1-2\varepsilon}{z^l}
	\frac{\Gamma(\alpha)}{\Gamma(-p+\varepsilon)\Gamma(\gamma)}
	F\left(\alpha,p+1-\varepsilon,\gamma;z^2\right).
\end{equation*}

Since $z<1$, the limit $R^<_{-1}(z)$ as $\varepsilon \to 0$ exists for all $p$.
If $p \ge 0$, then $R^<_{-1}(z)=0$ owing to the pole of $\Gamma(-p)$.
For $p<0$, the hypergeometric series reduces to a Jacobi polynomial of degree
$-p-1$:
\begin{equation}
	\label{eq:f_int}
	R^{\mathrm{reg},<}_{-1,k,l}(z) =
	\begin{cases}
		(2l-2p+1)\, z^l
		P^{(l+\frac12,1)}_{-p-1}(1-2z^2), & p < 0,\\
		0, & p \ge 0.
	\end{cases}
\end{equation}

For $z>1$, transformation~\eqref{eq:ext_int}, which is to be understood here
in the sense of the limit, gives
\begin{equation}
	\label{eq:f_ext}
	R^{\mathrm{reg},>}_{-1,k,l}(z) =
	\begin{cases}
		(2k+2p+3)\, z^{-k-4}
		P^{(k+\frac32,1)}_{p-1}(1-2z^{-2}), & p > 0,\\
		0, & p \le 0.
	\end{cases}
\end{equation}

The complete expression for the regular part of $R$ is given by
Eq.~\eqref{eq:reg} with $\lambda=-1$. For $p=0$ it vanishes for all
$z \ne 1$. This means that for $l=k+1$ the body forces are concentrated as a
single layer on the surface of the inclusion.

The expression $R_{-1}(1)$ must be understood in the sense of the limit
$\lambda \to -1+$; however, for $\lambda<0$ the hypergeometric series
in~\eqref{eq:hypergeo_int} diverges at $z=1$, and as $z \to 1$ it behaves as
$|1-z|^{\lambda}$~\cite{DLMF}.
% (see \cite[Eq.~15.4.23]{DLMF}).
Therefore, in addition to the regular part~\eqref{eq:reg}, $R_{-1}(z)$ also
contains the singular term
$R^{\mathrm{sing}}_{-1}(z)=(-1)^p\, \delta(z-1)$~\cite{Maximon,KR}.
To obtain it, one evaluates the action of $R_\lambda(z)$ on a test function as
$\lambda \to -1$, taking into account the above asymptotics of the
hypergeometric function. Alternatively, it can be obtained by
differentiating~\eqref{eq:reg} with respect to $z$ for $\lambda \to 0$, using
the well-known formula for derivatives of Bessel functions
(see~\cite[Eq.~10.6.5]{DLMF}).

Finally, we obtain
\begin{equation}
	\label{eq:f_tot}
	R_{-1,k,l} (z) = R^{\mathrm{reg},<}_{-1,k,l}(z)\,\theta(1-z) +
	R^{\mathrm{reg},>}_{-1,k,l}(z)\,\theta(z-1)
	+ (-1)^p\, \delta(z-1).
\end{equation}

\section{Conclusion}

Thus, for the cases considered, the values of the Weber--Schafheitlin integral
are expressed in terms of Jacobi polynomials and are given by the following
formulas:
\begin{itemize}
	\item for $\lambda=0$: Eqs.~\eqref{eq:reg}, \eqref{eq:e_int}, \eqref{eq:e_ext};
	\item for $\lambda=1$: Eqs.~\eqref{eq:reg}, \eqref{eq:u_int}, \eqref{eq:u_ext};
	\item for $\lambda=-1$: Eqs.~\eqref{eq:f_tot}, \eqref{eq:f_int}, \eqref{eq:f_ext}.
	In this case the integral is a distribution.
\end{itemize}

\textbf{Acknowledgments:} Parts of this manuscript were translated from Russian into English
using Claude Opus 5.5. The author reviewed and revised the translation for accuracy, terminology, and meaning, and takes full responsibility for the final text.

\end{document}